\documentclass[]{interact}

\usepackage{epstopdf}
\usepackage[caption=false]{subfig}

\usepackage[numbers,sort&compress]{natbib}
\bibpunct[, ]{[}{]}{,}{n}{,}{,}
\renewcommand\bibfont{\fontsize{10}{12}\selectfont}

\theoremstyle{plain}

\theoremstyle{definition}

\theoremstyle{remark}

\usepackage{amsthm}
\usepackage{amsmath, amssymb}
\usepackage{mathtools}
\mathtoolsset{showonlyrefs}
\usepackage{xcolor}

\def\ai#1{{\color{black}#1}} 
\def\gav#1{{\color{black}#1}} 
\def\dmiv#1{{\color{black}#1}} 
\def\mitya#1{{\color{black}#1}} 
\def\egor#1{{\color{black}#1}} 
\def\avg#1{{\color{black}#1}} 

\usepackage{hyperref}
\hypersetup{
    colorlinks=true,
    citecolor=blue,
    linkcolor=red, 
    filecolor=magenta,      
    urlcolor=cyan,
}

\usepackage{wrapfig}

\newcommand*{\R}{{\mathbb R}}
\newcommand\abs[1]{\left|#1\right|}

\usepackage{xcolor}
\usepackage{subfig}
\usepackage{hyperref}
\usepackage{textcomp}
\usepackage{nicefrac}

\usepackage{algorithm}
\usepackage{algpseudocode}
 \usepackage{algorithmicx}
\algnewcommand{\LeftComment}[1]{\Statex \(\triangleright\) #1}

\usepackage{xcolor}
\usepackage{bbm}

\DeclareMathOperator*{\argmin}{argmin}

\DeclareMathOperator*{\Argmin}{Argmin}

\newtheorem{thm}{Theorem}
\newtheorem{cor}{Corollary}
\newtheorem{prop}{\gav{Assumption}}
\newtheorem{rmk}{Remark}
\newtheorem{lem}{Lemma}

\begin{document}


\title{Adaptive Catalyst for Smooth Convex Optimization\footnote{The research of A. Gasnikov \gav{and V. Matyukhin} was supported by the Ministry of Science and Higher Education of the Russian Federation (Goszadaniye) \textnumero075-00337-20-03, project No. 0714-2020-0005.
} 
}

\author{
\name{Anastasiya Ivanova\textsuperscript{a,b} and Dmitry Pasechnyuk\textsuperscript{a} and Dmitry Grishchenko\textsuperscript{c} and Egor Shulgin\textsuperscript{a,d} and Alexander Gasnikov\textsuperscript{a,b,e} and Vladislav Matyukhin\textsuperscript{a}}
\affil{
\textsuperscript{a}Moscow Institute of Physics and Technology, Moscow, Russia;\\ 
\textsuperscript{b}National Research University Higher School of Economics, Moscow, Russia;\\
\textsuperscript{c}Universit\'e Grenoble Alpes, Grenoble, France;\\
\textsuperscript{d}King Abdullah University of Science and Technology, Thuwal, Saudi Arabia;\\
\textsuperscript{e}Institute for Information Transmission Problems, Moscow, Russia
}
}

\maketitle
\begin{abstract}
In this paper, we present a generic framework that allows accelerating almost arbitrary non-accelerated deterministic and randomized algorithms for smooth convex optimization
problems. The major approach of our \emph{envelope} is the same as in \textit{Catalyst} \citep{lin2015universal}: an accelerated proximal outer gradient method, which is used as an envelope for a non-accelerated inner method for the $\ell_2$ regularized auxiliary problem. Our algorithm has two key differences: $1)$ easily verifiable stopping criteria for inner algorithm; $2)$ the regularization parameter can be tuned along the way. As a result, the main contribution of our work is a new framework that applies to adaptive inner algorithms: Steepest Descent, Adaptive Coordinate Descent, Alternating Minimization. Moreover, in the non-adaptive case, our approach allows obtaining Catalyst without a logarithmic factor, which appears in the standard Catalyst \citep{lin2015universal,lin2018catalyst}.
\end{abstract}


\section{Introduction}
\normalsize
One of the main achievements in numerical methods for convex optimization is the development of accelerated methods \citep{nesterov2018lectures}. Until 2015 acceleration schemes for different convex optimization problems seem to be quite different to unify them. But starting from the work \citep{lin2015universal} in which universal acceleration technique (\textit{Catalyst}) was proposed, there appears a stream of subsequent works~\citep{lin2018catalyst,palaniappan2016stochastic,paquette2017catalyst,kulunchakov2019generic} that allows spreading Catalyst on monotone variational inequalities, non-convex problems, stochastic optimization problems. In all these works, the basic idea is to use an accelerated proximal algorithm as an outer envelope  \citep{rockafellar1976monotone} with non-accelerated algorithms for inner auxiliary problems. The main practical drawback of this approach is the requirement to choose a regularization parameter such that the conditional number of the auxiliary problem becomes became $O(1)$. To do that, we need to know the smoothness parameters of the target that are not typically free available. 

An alternative accelerated proximal envelop\dmiv{e} \citep{parikh2014proximal} was proposed in the paper \citep{monteiro2013accelerated}. The main difference with standard accelerated proximal envelops is the adaptability of the scheme \citep{monteiro2013accelerated}. Note, that this scheme allows also to build (near) optimal tensor (high-order) accelerated methods \citep{gasnikov2017universal,nesterov2018lectures,gasnikov2019near,gasnikov2019reachability,wilson2019accelerating}. That is, the ``acceleration'' potential of this scheme seems to be the best known for us for the moment. So the main and rather simple idea of this paper can be formulated briefly as follows: \textbf{To develop adaptive Catalyst, we replace the accelerated proximal envelope with a fixed regularization parameter \citep{parikh2014proximal,lin2018catalyst} on the adaptive accelerated proximal envelope from \citep{monteiro2013accelerated}}. 

\gav{In} Section~\ref{2} \gav{, we describe adaptive Catalyst envelop\dmiv{e} -- Algorithm~\ref{alg:adaptive_catalyst} and generalized Monteiro-Svaiter theorem from \cite{monteiro2013accelerated} to set out how to do this envelope work. We emphasize that the proof of the theorem contains as a byproduct the new theoretical analysis of the stopping criterion for the inner algorithm \eqref{inner_est}. This stopping criterion allows one to show that the proposed envelope in non-adaptive mode is log-times  better (see Corollary~\ref{cor1}) in the total number of oracle calls of the inner method (we'll measure the complexity of the envelope in such sense) in comparison with all other envelopes known for us.}

By using this adaptive accelerated proximal envelope,  we propose in Section~\ref{3} an accelerated variant of steepest descent \citep{polyak1987introduction,gasnikov2017universal} as an alternative to A. Nemirovski accelerated steepest descent (see \citep{nesterov2018primal,diakonikolas2019conjugate} and references therein), adaptive accelerated variants of alternating minimization procedures \citep{beck2017first} as an alternative to \citep{diakonikolas2018alternating,guminov2019accelerated,tupitsa2019alternating} and adaptive accelerated coordinate descent \citep{nesterov2012efficiency}. For the last example, as far as we know, there were no previously complete adaptive accelerated coordinate descent. The most advanced result in this direction is the work \citep{fercoq2015accelerated} that applies only to the problems with increasing smoothness parameter along the iteration process. For example, for the target function like $f(x) = x^4$, this scheme doesn't recognize that smoothness parameters (in particular Lipschitz gradient constant) tend to zero along the iteration process.

In Section~\ref{4} we describe numerical experiments with the steepest descent, adaptive coordinate descent, \gav{alternating minimization and local SGD}. \gav{We try to accelerate these methods by the envelope, described in Algorithm~\ref{alg:adaptive_catalyst}.}

We hope that the proposed approach allows accelerating not only adaptive on their own procedures, but also many other different non-accelerated non-adaptive randomized schemes by settings on general smoothness parameters of target function that can be difficult to analyze patently \citep{gower2019sgd,gazagnadou2019optimal,gorbunov2019unified}.

\gav{The first draft of this paper appeared in arXiv in November 2019. Since that time, this paper has developed (and cited) in different aspects. The main direction is a convenient (from the practical \eqref{MS_cond} and theoretical \eqref{inner_est} point of view) criterion of stopping the inner algorithm that is wrapped in  accelerated proximal envelope. We emphasize that our contribution in this part is not a new accelerated proximal envelop\dmiv{e} (we use the well-known envelope \cite{monteiro2013accelerated}), but we indicate that this envelope is better than the other ones due to the new theoretical analysis of its inner stopping criteria that lead us from \eqref{MS_cond} to \eqref{inner_est}. Although this calculation looks simple enough, to the best of our knowledge, this was the first time when it was provably developed an accelerated proximal envelope that required to solve the auxiliary problem with prescribed relative accuracy in argument $\simeq 1/5$. Since the auxiliary problem is smooth and strongly convex, this observation eliminates the logarithmic factor (in the desired accuracy) in the complexity estimate for such an envelope in comparison with all known analogues. Note that this small observation will have a remarkable influence on the development of accelerated algorithms. A close stopping condition, for instance, arises in the following papers \cite{doikov2020contracting,doikov2020inexact} that developed (sub-)optimal accelerated tensor method based on accelerated proximal envelopes. The proposed ``logarithm-free'' envelope allows one to improve the best known bounds  \cite{lin2020gradient} for strongly convex-concave saddle-point problems (with different constants of strong convexity and concavity) on logarithmic factor \cite{dvinskikh2020accelerated_}. Composite variant of this envelop\dmiv{e} also allows one to develop ``logarithm-free'' gradient sliding-type methods\footnote{Note, that \cite{dvinskikh2020accelerated} contains variance reduction  \cite{shalev2014accelerated,allen2016optimal} generalization (with non proximal-friendly composite) of proposed in this paper scheme.}  \cite{ivanova2020oracle,dvinskikh2020accelerated,dvinskikh2020accelerated_} and its tensor generalizations \cite{dvurechensky2020optimal}. Moreover, one of the variants of the hyper-fast second-order method was also developed based on this envelope \cite{kamzolov2020near}. Though this envelope had known before, it seems that the original idea of this paper to use this envelope in Catalyst type procedures and new (important from the theoretical point of view) reformulation of stopping criteria for inner algorithm \eqref{inner_est} has generated a large number of applications, some of them mentioned in this paper, the others can be found in the literature cited above. As an important example, we show in section~\ref{RACDM} that the developed envelop\dmiv{e} with non-accelerated coordinate descent method for auxiliary problem works much better in theory (and better in practice) than all known direct accelerated coordinate-descent algorithms for sparse soft-max type problem. Before this article, this was an open problem, how to beat standard accelerated coordinate-descent algorithms that don't allow one to take into account sparsity of the problem for soft-max type functional \cite{GasnikovBook,pasechnyuk2021accelerated}.}

\gav{The other contribution of this paper is an adaptive choice of the smooth parameter. Since our approach requires two inputs (lower bound $L_d$ and upper bound $L_u$ for the unknown smoothness parameter $L$), it's hardly possible to call it ``adaptive''. Moreover, the greater the discrepancy between these two parameters, the worthier is our adaptive envelope in theory. But almost all of our experiments demonstrate low sensitivity to these parameters rather than to real smoothness parameter. But even for such a ``logarithm-free'' and adaptive envelope, we expect that typically the direct adaptive accelerated procedures will work better than its Catalyst type analogues. It was recently demonstrated in the following work \cite{tupitsa2020accelerated} for accelerated alternating minimization procedure. But even to date, there are problems in which one can expect that firstly optimal accelerated algorithms will be developed by using Catalyst type procedures rather than direct acceleration. Recent advances in saddle-point problems \cite{lin2020gradient,dvinskikh2020accelerated_,yang2020catalyst} and decentralized distributed optimization\footnote{Note, that the results of these papers were further reopened by using direct acceleration \cite{kovalev2020optimal,li2020optimal}.} \cite{li2020revisiting,hendrikx2020dual} confirm this thought. We expect that for Homogeneous Federated Learning architectures, Accelerated Local SGD can be developed (see \cite{woodworth2020local} for the state of the art approach) by using Catalyst-type envelop\dmiv{e} with SCAFFOLD version of local SGD algorithm \cite{karimireddy2019scaffold}. As far as we know, it's still an open problem to build Accelerated local SGD as an inner algorithm. In this paper, we demonstrate some optimistic experiments in this direction.}

\section{The Main Scheme}\label{2}
Let us consider the following minimization problem 
\begin{equation}\label{init_problem}
    \min\limits_{y \in \gav{\R^n}} f(y),
\end{equation}
where $f(y)$ is a convex function, and its gradient is Lipschitz continuous w.r.t. $\|\cdot\|_{2}$ with the constant $L_{f}$:
\begin{eqnarray*}
     \| \nabla f(x) - \nabla f(y)\|_{2} \leq L_{f} \|x - y \|_{2}.
\end{eqnarray*} 
\gav{We denote $x_{\star}$ a solution of \eqref{init_problem}.}

To propose the main scheme of the algorithm we need to define the following functions:
\begin{eqnarray*}
    F_{L,x}(y)& = & f(y) + \tfrac{L}{2}||y-x||^2_2,\\
    f_{L}(x)& = &\min\limits_{y \in \gav{\R^n}}F_{L,x}(y) = F_{L,x}(y_{L}(x)),
\end{eqnarray*}
then the function $F_{L,x}(y)$ is $L$--strongly convex, and its gradient is Lipschitz continuous w.r.t. $\|\cdot\|_{2}$ with the constant $(L+L_{f})$. So, the following inequality holds 
\begin{equation}
\label{lip_1}
    || \nabla F_{L,x}(y_2) - \nabla F_{L,x}(y_1)||_2 \leq  (L+L_{f})||y_1 - y_2||_2.
\end{equation}
Due to this definition, for all $L \geq 0$ we have that $f_{L}(x) \leq f(x)$ and the convex function $f_{L}(x)$ has a Lipschitz-continuous gradient with the Lipschitz constant $L$. Moreover, according to \citep{polyak1987introduction} [Theorem 5, ch. 6], since 
\begin{eqnarray*}
    x_\star \gav{\in} 
    \Argmin\limits_{x \in \gav{\R^n}}f_{L}(x),
\end{eqnarray*}
we obtain 
\begin{eqnarray*}
    x_\star \in \Argmin\limits_{x \in \gav{\R^n}}f(x) \quad \text{and } f_L(x_\star)=f(x_\star).
\end{eqnarray*}
Thus, instead of the initial problem~\eqref{init_problem}, we can consider the Moreau--Yosida regularized problem
\begin{equation}\label{new_problem}
    \min\limits_{x \in \gav{\R^n}} f_L(x).
\end{equation}
Note that the problem~\eqref{new_problem} is an ordinary problem of smooth convex optimization. 
Then the complexity of solving the problem~\eqref{new_problem} up to the accuracy $\varepsilon$ with respect to the function 
using the Fast Gradient Method (FGM) \citep{nesterov2018lectures} can be estimated as follows $O\left( \sqrt{\tfrac{LR^2}{\varepsilon}} \right)$. The `complexity' means here the number of oracle calls. Each oracle call means calculation of $\nabla f_L(x)= L(x - y_L(x))$, where $y_L(x)$ is the exact solution of the auxiliary problem $\min\limits_{y \in \gav{\R^n}}F_{L,x}(y)$.

Note that the smaller the value of the parameter $L$ we choose, the smaller is the number of oracle calls (outer iterations). However, at the same time, this increases the complexity of solving the auxiliary problem at each iteration. 

At the end of this brief introduction to standard accelerated proximal point methods, let us describe the step of ordinary (proximal) gradient descent (for more details see \citep{parikh2014proximal})
\begin{equation*}
    x^{k+1} = x^{k} - \tfrac{1}{L}\nabla f_L(x^k) = x^{k} - \tfrac{L}{L}(x^k - y_L(x^k))= y_L(x^k).
\end{equation*}
To develop an adaptive proximal accelerated envelop, we should replace standard FGM \citep{nesterov2018lectures} on the following adaptive variant of FGM Algorithm~\ref{alg:adaptive_catalyst}, introduced by \citep{monteiro2013accelerated}
for smooth convex optimization problems.

\begin{algorithm}[tb]
\caption{Monteiro--Svaiter algorithm} 
\label{alg:adaptive_catalyst}
\begin{algorithmic}
	\State {\bf Parameters:} 
	$z^0, y^0, A_0 = 0$ 
	\For{$k = 0,1, \ldots, N-1$}
	\State Choose $L_{k+1}$ and $y^{k+1}$ such that
	\begin{eqnarray*}
	\|\nabla F_{L_{k+1}, x^{k+1}} (y^{k+1})\|_2 \leq \tfrac{L_{k+1}}{2} \|y^{k+1} - x^{k+1}\|_2,
	\end{eqnarray*} 
	
	where
	\begin{eqnarray*}
	    a_{k+1} &=& \tfrac{1/L_{k+1} + \sqrt{1/L^2_{k+1} + 4 A_k/L_{k+1}}}{2},\\
	    A_{k+1} &=& A_k + a_{k+1},\\
	    x^{k+1}& =& \tfrac{A_{k}}{A_{k+1}} y^{k} + \tfrac{a_{k+1}}{A_{k+1}} z^{k}
	\end{eqnarray*}
	
	\State $z^{k+1}=z^{k}-a_{k+1} \nabla f\left(y^{k+1}\right)$
	\EndFor

\State {\bf Output:} $y^{N}$
\end{algorithmic}	
\end{algorithm}
The analysis of the algorithm is based on the following theorem.
\begin{thm}[{Theorem 3.6~\citep{monteiro2013accelerated}}]
\label{Th:M-S-conv}
Let sequence $(x^k,y^k,z^k)$, $k \geq 0$ be generated by Algorithm \ref{alg:adaptive_catalyst} and define $R:=\left\| {y^0-x_{\star} } \right\|_2$. Then, for all $N \geq 0$, 
\begin{align}
\label{eq11}
\tfrac{1}{2}\left\|{z^N-x_{\star} }\right\|_2^2 &+A_N \cdot \left( {f\left( 
{y^N} \right)-f\left( {x_{\star} } \right)} 
\right) +\tfrac{1}{4}\sum\limits_{k=1}^N {A_k L_{k} \left\| {y^k-x^{k}} 
\right\|_2^2 } \le \tfrac{R^2}{2},
\end{align}
\begin{equation}
\label{eq12}
f\left( {y^N} \right)-f\left( {x^{\star} } \right)\le \tfrac{R^2}{2A_N },
\quad
\left\| {z^N-x_{\star} } \right\|_2 \le R,
\end{equation}
\begin{equation}
\label{eq13}
\sum\limits_{k=1}^N {A_k L_{k} \left\| {y^k-x^{k}} \right\|_2^2 } \le 
2R^2.
\end{equation}
\end{thm}

We also need the following Lemma.
\begin{lem}[{Lemma 3.7a~\citep{monteiro2013accelerated}}] \label{pr:prop2}
    Let sequences $\left\{ {A_k ,L_k } \right\}$, $k\geq 0$ be generated by Algorithm \ref{alg:adaptive_catalyst}. Then, for all $N\geq 0$,
\begin{equation}
\label{eq14}
A_N \ge \tfrac{1}{4}\left( {\sum\limits_{k=1}^N {\tfrac{1}{\sqrt {L_{k} } }} 
} \right)^2.
\end{equation}    
\end{lem}

Let us define non-accelerated method $\mathcal{M}$ that we will use to solve auxiliary problem. 

\begin{prop}
\label{main_prop}
The convergence rate (after $t$ iterations / oracle calls) for the method $\mathcal{M}$ for problem
$$\min\limits_{y \in \R^n}F(y)$$
can be written in the general form as follows: with probability at least $ 1-\delta$ holds (for randomized algorithms, like Algorithm~\ref{algo:CD}, this estimates holds true with high probability)\footnote{For deterministic algorithms we can skip ``with probability at least $1 - \delta$'' and factor ``$\ln\tfrac{N}{\delta}$''.} 
\begin{equation*}
    F(y_t) - F(y_\star) \,  =\,O\left(L_{F} R\ai{_y}^2 \ln\tfrac{t}{\delta}\right)
    \min \left\{\tfrac{C_n}{t}, \exp{\left(-\tfrac{\mu_{F} t}{C_n L_{F}} \right)} \right\},
\end{equation*}
where $y_\star$ is the solution of the problem, $R_{y} = ||y^0 - y_\star||_2$, function $F$ is $\mu_{F}$--strongly convex and $L_{F}$ is a constant which characterized smoothness of function $F$. 
\end{prop}

Typically $C_n=O(1)$  for the standard full gradient first order methods, $C_n=O(p)$, where $p$ is a number of blocks, for alternating minimization with $p$ blocks and $C_n=O(n)$ for gradient free or coordinate descent methods, where $n$ is dimension of $y$. See the references in next Remark for details.

\begin{rmk} Let us clarify what we mean by a constant $L_{F}$ which characterized smoothness of function $F$. Typically for the first order methods this is just the Lipschitz constant of gradient $F$ (see, \citep{polyak1987introduction,de2017worst} for the steepest descent and  \citep{karimi2016linear,diakonikolas2018alternating,tupitsa2019alternating} for alternating minimization); for gradient free methods like Algorithm~\ref{algo:CD} this constant is the average value of the directional smoothness parameters, for gradient free methods see \citep{duchi2015optimal,gasnikov2016gradient-free,shamir2017optimal,bayandina2017gradient-free,dvurechensky2018accelerated,dvurechensky2018accelerated2}, for coordinate descent methods see  \citep{nesterov2012efficiency, wright2015coordinate,nesterov2017efficiency} and for more general situations see \citep{gower2019sgd}.
\end{rmk}

\begin{rmk} Note that in \gav{Assumption}~\ref{main_prop} the first estimate corresponds to the estimate of the convergence rate of the method  $\mathcal{M}$ for convex problems. And the second estimate corresponds to the estimate for strongly convex problems.
\end{rmk}

Our main goal is to propose a scheme to accelerate methods of this type. But note that we apply our scheme only to degenerate convex problems since it does not take into account the strong convexity of the original problem.

Denote $F^{k+1}_{L, x} (\cdot) \equiv F_{L_{k+1}, x^{k+1}} (\cdot)$. Based on Monteiro--Svaiter accelerated proximal method we propose Algorithm~\ref{alg1}. 

\begin{algorithm}[tb]
\caption{Adaptive Catalyst}
\label{alg1}
\begin{algorithmic}
\State {\bf Parameters:}  Starting point $x^0=y^0=z^0$;
initial guess $L_0 > 0$; parameters $\alpha > \beta \gav{\gtrsim} \gamma > \gav{1}$; optimization method $\mathcal{M}$, $A_0=0$.
\For{$k = 0, 1, \ldots,N-1$}
    \State{$L_{k+1} =  \beta \cdot\min\left\{ \alpha L_{k},L_u\right\}$}
     \State{$r = 0$}
    \Repeat 
    \State{$r: = r+1$}
     \State{$L_{k+1}:=  \max\left\{L_{k+1}/ \beta,L_d\right\}$}
        \State Compute 
           \begin{eqnarray*}
	    a_{k+1} &=& \tfrac{1/L_{k+1} + \sqrt{1/L^2_{k+1} + 4 A_k/L_{k+1}}}{2},\\
	    A_{k+1} &=& A_k + a_{k+1},\\
	    x^{k+1}& =& \tfrac{A_{k}}{A_{k+1}} y^{k} + \tfrac{a_{k+1}}{A_{k+1}} z^{k}.
	      \end{eqnarray*}
        \State Compute an approximate solution of the following problem with auxiliary non-accelerated method $\mathcal{M}$
	$$ y^{k+1} \approx  \argmin_{y}  F^{k+1}_{L, x} (y)  \, :$$ 
    \ai{ By running $\mathcal{M}$ with starting point $x^{k+1}$ and output point $y^{k+1}$ we wait  $N_r$ iterations to fulfill adaptive stopping criteria 
	$$ \|\nabla F^{k+1}_{L, x}  (y^{k+1})\|_2 \leq \tfrac{L_{k+1}}{2} \|y^{k+1} - x^{k+1}\|_2. $$} 
    \Until{$r>1$ and $N_{r} \ge \gamma \cdot N_{r-1}$ or $L_{k+1} = L_d$}
    \State $z^{k+1}=z^{k}-a_{k+1} \nabla f\left(y^{k+1}\right)$
    \EndFor
\State {\bf Output:} $y^{N}$
\end{algorithmic}
\end{algorithm}

Now let us prove the main theorem about the convergence rate of the proposed scheme. \ai{Taking into account that $\tilde{O}(\cdot)$ means the same as $O(\cdot)$ up to a logarithmic factor, based on the Monteiro--Svaiter Theorem~\ref{Th:M-S-conv} we can introduce the following theorem:}

\begin{thm}
\label{main_theo}
Consider Algorithm~\ref{alg1} with $0 < L_d < L_u$ for solving problem~\eqref{init_problem}, where $Q = \R^n$, with auxiliary (inner) non-accelerated algorithm (method) $\mathcal{M}$ that satisfy \gav{Assumption}~\ref{main_prop} with constants $C_n$ and $L_f$ such that $L_d\le L_f \le L_u$.

 Then the total complexity\footnote{The number of oracle calls (iterations) of auxiliary method $\mathcal{M}$ that required to find $\varepsilon$ solution of \eqref{init_problem} in terms of functions value.} of the proposed Algorithm~\ref{alg1} with inner method $\mathcal{M}$ is $$\tilde{O}\left(C_n \cdot \max\left\{\sqrt{\tfrac{L_u}{L_f}},\sqrt{\tfrac{L_f}{L_d}}\right\}\cdot \sqrt{\tfrac{L_f R^2}{\varepsilon}}\right)$$
with probability at least $1 - \delta$.
\end{thm}
\begin{proof}

Note that the Monteiro--Svaiter (M-S) condition
\begin{equation}\label{MS_cond}
    \|\nabla F^{k+1}_{L, x}  (y^{k+1})\|_2 \leq \tfrac{L_{k+1}}{2} \|y^{k+1} - x^{k+1}\|_2 
\end{equation}
instead of the exact solution $y^{k+1}_\star = y_{L_{k + 1}} (x^{k+1})$ of the auxiliary problem, for which
\begin{eqnarray*}
    \|\nabla F^{k+1}_{L, x}  (y^{k+1}_\star)\|_2 = 0, 
\end{eqnarray*}
allows to search inexact solution that satisfies the condition~\eqref{MS_cond}. 

Since $y^{k+1}_\star$ \gav{is} the solution of the problem  $\min\limits_{y}  F^{k+1}_{L, x} (y)$, the $\nabla F^{k+1}_{L, x}  (y^{k+1}_\star) = 0$. Then, using inequality~\eqref{lip_1} we obtain 
\begin{equation}
\label{eq:1}
    || \nabla F^{k+1}_{L, x} (y^{k+1})||_2 \leq  (L_{k+1}+L_{f})||y^{k+1} - y^{k+1}_\star||_2.
\end{equation}
Using the triangle inequality we have 
\begin{equation}
  \label{eq:2}
    ||x^{k+1} - y^{k+1}_\star||_2 - ||y^{k+1} -  y^{k+1}_\star||_2 \leq || y^{k+1} - x^{k+1} ||_2.  
\end{equation}

Since r.h.s. of the inequality~\eqref{eq:2} coincide with the r.h.s. of the M-S condition and l.h.s. of the inequality~\eqref{eq:1}
coincide with the l.h.s. of the M-S condition up to a multiplicative factor $L_{k+1}/2$, one can conclude that if the inequality 
\begin{equation}
\label{inner_est}
    ||y^{k+1} - y^{k+1}_\star||_2 \leq \tfrac{L_{k+1}}{3 L_{k+1} + 2L_{f}}  ||x^{k+1} - y^{k+1}_\star||_2
    \end{equation}
holds, the M-S condition holds too. 

To solve the auxiliary problem $\min\limits_{y}  F_{L_{k+1}, x^{k+1}}(y)$ we use non-accelerated method $\mathcal{M}$. Using \gav{Assumption}~\ref{main_prop} with probability $\geq 1 - \tfrac{\delta}{N}$ (where $N$ is the total number of the Catalyst's steps), we obtain that the convergence rate (after $t$ iterations of $\mathcal{M}$, see \gav{Assumption}~\ref{main_prop})

\begin{equation*}
    F^{k+1}_{L, x} (y^{k+1}_t) - F^{k+1}_{L, x} (y^{k+1}_\star) \notag =
    O\left((L_{f}+L_{k+1}) R_{k+1}^2 \ln\tfrac{N t}{\delta}\right)
     \exp{\left(-\tfrac{L_{k+1} t}{C_n(L_{f}+L_{k+1})} \right)}.
\end{equation*}
Note, that $R_{k+1} = ||x^{k+1} - y^{k+1}_\star||_2$ since $x^{k+1}$ is a starting point. 

Since $F^{k+1}_{L, x} (\cdot)$ is $L_{k+1}$-strongly convex function, the following inequality holds \citep{nesterov2018lectures}
\begin{eqnarray*}
 \tfrac{L_{k+1}}{2}||y^{k+1}_t- y^{k+1}_\star||^2_2 
  \leq F^{k+1}_{L, x} (y^{k+1}_t) - F^{k+1}_{L, x} (y^{k+1}_\star).
\end{eqnarray*}
Thus,
\begin{equation}\label{y_acc}
||y^{k+1}_t - y^{k+1}_\star||_2 \leq  O\left( \sqrt{\tfrac{(L_{f}+L_{k+1}) R_{k+1}^2 }{L_{k+1}}\ln\tfrac{N t}{\delta}}\right) \exp{\left(-\tfrac{L_{k+1}t }{2C_n(L_{f}+L_{k+1})} \right)}.
\end{equation}
From \gav{\eqref{inner_est}, \eqref{y_acc} and the fact that we start $\mathcal{M}$ at $x^{k+1}$}, we obtain that the complexity \gav{$T$ (number of iterations of $\mathcal{M}$)} of solving the auxiliary problem with probability at least $1 - \tfrac{\delta}{N}$ is \gav{determined from}
\gav{
\begin{equation}
\label{ForCor1}
O\left(R_{k+1}\sqrt{\tfrac{(L_{f}+L_{k+1})}{L_{k+1}}\ln\tfrac{N T}{\delta}}\right) \exp\left(-\tfrac{L_{k+1}T }{2C_n(L_{f}+L_{k+1})} \right) 
\simeq \tfrac{L_{k+1}}{3 L_{k+1} + 2L_{f}}  R_{k+1},
\end{equation}
}
\gav{hence}
\begin{equation}
\label{comp_aux}
    T = \tilde{O}\left(C_n\tfrac{(L_{k+1}+ L_{f}) }{L_{k+1}}\right).
\end{equation}
Since \gav{we use in \eqref{comp_aux} $\tilde{O}(\text{ })$ notation,} we can consider $T$ to be the estimate that \gav{corresponds to the} total complexity of auxiliary problem including all inner restarts on $L_{k+1}$.

  Substituting inequality~\eqref{eq14} into estimation \eqref{eq13} we obtain 
 \begin{eqnarray*}
     f(y^N) - f(x_\star) \leq \tfrac{2R^2}{\left(\sum\limits_{k=1}^N \tfrac{1}{\sqrt{L_{k}}}\right)^2}.
 \end{eqnarray*}
Since the complexity of the auxiliary problem with probability at least $1 - \tfrac{\delta}{N}$ \ai{is} $T$ we assume that in the worst case all $L_{k+1}$ are equal. Then the worst case we can estimate as the following optimization problem
$$
    \max\limits_{L_d \leq L \leq L_u} \tfrac{L + L_f}{L} \sqrt{\tfrac{LR^2}{\varepsilon}},
$$
Obviously, the maximum is reached at the border.
 So, using union bounds inequality over all $N$ iterations of the Catalyst we can estimate the complexity in the worst two cases as follows:
\begin{itemize}
    \item If all $L_{k+1} = L_{d} \leq  L_{f}$ (at each iteration we estimate the regularization parameter as lower bound), then $\tfrac{(L_{k+1}+ L_{f} ) }{L_{k+1}} \approx \tfrac{ L_{f}  }{L_{k+1}}$ and total complexity with  probability $\geq 1 - \delta$ is 
    $$
\tilde{O}\left(C_n \tfrac{ L_{f}  }{L_{d}} \sqrt{ \tfrac{L_{d}R^2}{\varepsilon}} \right) = \tilde{O}\left(C_n \sqrt{ \tfrac{L_{f}}{L_d}} \cdot \sqrt{ \tfrac{L_{f}R^2}{\varepsilon}}\right).
$$ 
\item If all $L_{k+1} = L_{u} \geq  L_{f}$ (at each iteration we estimate the regularization parameter as upper bound), then $\tfrac{(L_{k+1}+ L_{f} ) }{L_{k+1}} \approx 1$ and total complexity with  probability $\geq 1 - \delta$ is 
    $$
\tilde{O}\left(C_n  \sqrt{ \tfrac{L_{u}R^2}{\varepsilon}} \right) = \tilde{O}\left(C_n \sqrt{ \tfrac{L_{u}}{L_f}} \cdot \sqrt{ \tfrac{L_{f}R^2}{\varepsilon}}\right).
$$ 
\end{itemize}
Then, using these two estimations we obtain the result of the theorem.
\end{proof}

Note that this result shows that such a procedure will works not worse than standard Catalyst \citep{lin2015universal,lin2018catalyst} up to a factor $\tilde{O}\left(\max\left\{\sqrt{\frac{L_u}{L_{f}}},\sqrt{\frac{L_{f}}{L_d}}\right\}\right)$ independent on the stopping criteria in the restarts on $L_{k+1}$.

Since the complexity of solving the auxiliary problem is proportional to $\frac{(L_{k+1}+ L_{f} )C_n }{L_{k+1}}$,
when we reduce the parameter $L_{k+1}$ so that $L_{k+1} < L_{f}$ the complexity of solving an auxiliary problem became growth exponentially. Therefore, as the stopping criterion of the inner method, we select the number of iterations $N_{t}$ compared to the number of iterations $N_{t-1}$ at the previous restart $t-1$. This means that if  $N_{t} \leq \gamma N_{t-1}$ then the complexity begins to grow exponentially and it is necessary to go to the next iteration of the external method. By using such adaptive rule we try to recognize the best possible value of $L_{k+1}\simeq L_f$. The last facts are basis of standard Catalyst approach \citep{lin2015universal,lin2018catalyst} and have very simple explanation. To minimize the total complexity we should take parameter $L_{k+1}\equiv L$  such that $$\min_{L}\sqrt{\tfrac{LR^2}{\varepsilon}}\cdot\tilde{O}\left(\tfrac{L_f + L}{L}\right). $$ This leads us to $L_{k+1}\simeq L_f$. 

Note that also in non-adaptive case (if we choose all $L_{k+1} \equiv L_{f}$) we can obtain the following corollary from the Theorem~\ref{main_theo}.
\begin{cor}\label{cor1}
If we consider Algorithm~\ref{alg1} with $L_{k+1} \equiv L_{f}$ for solving problem~\eqref{init_problem}, then the total complexity of the proposed Algorithm~\ref{alg1} with inner \ai{non-randomized} method $\mathcal{M}$ is  
\begin{equation}
\label{cor_1_est}
    O\left(C_n \sqrt{\tfrac{L_f R^2}{\varepsilon}}\right).
\end{equation}
\end{cor}
\begin{proof}
Using \gav{\eqref{ForCor1} without $\ln\tfrac{N T}{\delta}$ factor (since $\mathcal{M}$ is non-randomized)}
\gav{we derive}
that the complexity of the auxiliary problem \gav{is (see also~\eqref{comp_aux})}
\begin{equation*}
    T = O\left(C_n\tfrac{(L_{k+1}+ L_{f}) }{L_{k+1}} \cdot \ln\left(\tfrac{3L_{k+1} + 2L_{f}}{L_{k+1}} \right) \right)
\end{equation*}
And since \gav{we choose  $L_{k+1} \equiv L_{f}$,}
\gav{$$\tfrac{3L_{k+1} + 2L_{f}}{L_{k+1}} = 5$$}
Then the complexity of the auxiliary problem is $ T = O\left(\gav{C_n
}\right)$. Using this estimate\gav{,} 
we obtain that the total complexity is~\eqref{cor_1_est}.
\end{proof}

If method $\mathcal{M}$ is randomized we have the additional factor $\ln\tfrac{N T}{\delta}\simeq \ln\left(\tfrac{1}{\delta\varepsilon}\right)$. Hence, \eqref{cor_1_est} changes: with probability at least $1-\delta$
\begin{equation*}
    O\left(C_n\ln\left(\tfrac{1}{\delta\varepsilon}\right) \sqrt{\tfrac{L_f R^2}{\varepsilon}}\right).
\end{equation*}

\ai{Note that in the standard Catalyst approach~\citep{lin2015universal,lin2018catalyst} the total complexity is $O\left(C_n\ln\left(\tfrac{1}{\delta\varepsilon}\right)\sqrt{\tfrac{L_{f}R^2}{\varepsilon}}\cdot \ln\left(\tfrac{1}{\varepsilon'}\right)\right)$,} where $\varepsilon'=\text{Poly}(\varepsilon)$ is the relative accuracy of solving the auxiliary problem at each iteration. From this we get that choosing the stopping criterion for the inner method as the criterion from the Algorithm~\ref{alg1} we can get the Catalyst without a logarithmic cost $\ln\left(\tfrac{1}{\varepsilon'}\right)$.
 It seems that such  variant of Catalyst can be useful in many applications. For example, as universal envelope for non accelerated asynchronized centralized distributed algorithms \citep{mishchenko2018delay}.

\section{Applications}\label{3}
In this section, we present a few examples of algorithms that we consider as inner solvers. \gav{Most} of them have an adaptive structure. \gav{It's natural to apply adaptive envelop\dmiv{e} to adaptive algorithms since the developed methods keep adaptability.} 
\subsection{Steepest Descent}

Consider the following problem 
\begin{equation*}
\min_{x\in \R^n}  f(x),
\label{eq:PrSt}
\end{equation*}
where $f(x)$ is a $L_f$--smooth convex function (its gradient is Lipschitz continuous w.r.t. $\|\cdot\|_{2}$ with the constant $L_{f}$).

To solve this problem, let us consider the general gradient descent update rule
\begin{equation*}\label{eq:gd}
     x^{k+1} = x^k - h_k\nabla f(x^k).
\end{equation*}
In \citep{polyak1987introduction} it was proposed an adaptive way to select $h_k$ as following  (see also \citep{de2017worst} for precise rates of convergence)
\begin{equation*}\label{eq:Steep}
    h_k = \argmin_{h\in\mathbb{R}} f(x^k - h\nabla f(x^k)).
\end{equation*}

\begin{algorithm}[ht]
\caption{Steepest descent}
\label{algo:SD}
\begin{algorithmic}
   \State {\bf Parameters:} Starting point $x^0$.
    \For{$k = 0,1,\ldots,N-1$}
	\State Choose $h_k = \argmin_{h\in\mathbb{R}}f(x^k - h\nabla f(x^k))$
	\State Set $x^{k+1} = x^k - h_k\nabla f(x^k)$
    \EndFor
\State {\bf Output:} $x^N$
\end{algorithmic}
\end{algorithm}

In contrast with the standard selection $h_k\equiv\frac{1}{L_f}$ for $L$-smooth functions $f$, in this method there is no need to know smoothness constant of the function. It allows to use this method for the smooth functions $f$ when $L_f$ is unknown (or expensive to compute) or when the global $L_f$ is much bigger than the local ones along the trajectory. 

On the other hand, as far as we concern, there is no direct acceleration of the steepest descent algorithm. Moreover, it is hard to use Catalyst with it as far as acceleration happens if $L_k$ ($\kappa$ in Catalyst article notations) is selected with respect to $L_f$ and the scheme does not support adaptivity out of the box. Even if global $L_f$ is known, the local smoothness constant could be significantly different from it that will lead to the worse speed of convergence.

Note that for Algorithm~\ref{algo:SD} the \gav{Assumption}~\ref{main_prop} holds with $C_n = O(1)$ and $L_{f}$ is the Lipschitz constant of the gradient of function $f$.

\subsection{Random Adaptive Coordinate Descent Method}\label{RACDM}
Consider the following unconstrained problem 
\begin{equation*}
\min_{x\in \R^n}  f(x).
\end{equation*}
Now we assume directional smoothness for $f$, that is there exists $\beta_1, \dots, \beta_n$ such that for any $ x \in \mathbb{R}^n, u \in \mathbb{R}$
$$
\left|\nabla_{i} f\left(x+u e_{i}\right)-\nabla_{i} f(x)\right| \leq  \beta_{i}|u|, \quad  i =1, \ldots, n,
$$
where $\nabla_{i} f(x) = \partial f(x) / \partial x_i.$
For twice differentiable $f$ it equals to $(\nabla^{2} f(x))_{i, i} \leq \beta_{i}.$
Due to the fact that we consider the situation when smoothness constants are not known, we use such a dynamic adjustment scheme from~\citep{nesterov2012efficiency,wright2015coordinate}.

\begin{algorithm}[ht]
\caption{RACDM}
\label{algo:CD}
\begin{algorithmic}
   \State {\bf Parameters:} Starting point $x^0$;\\
   lower bounds $\hat{\beta}_{i}:=\beta_{i}^{0} \in\left(0, \beta_{i}\right], i=1, \dots, n$
    \For{$k = 0,1,\ldots,N-1$}
	\State Sample $i_k \sim \mathcal{U}[1, \dots, n]$ 
	\State Set $x^{k+1} = x^k - \hat{\beta}_{i_{k}}^{-1} \cdot \nabla_{i_{k}} f\left(x^{k}\right) \cdot e_{i_{k}}$
	\State \textbf{While} $\nabla_{i_{k}} f(x^{k}) \cdot \nabla_{i_{k}} f(x^{k+1}) < 0$\, \textbf{do}
	$$\left\{ \hat{\beta}_{i_{k}}=2 \hat{\beta}_{i_{k}}, \quad x^{k+1}=x^{k}-\hat{\beta}_{i_{k}}^{-1} \cdot \nabla_{i_{k}} f\left(x^{k}\right) \cdot e_{i_{k}} \right\}$$
	\State Set $\beta_{i_{k}} = \frac{1}{2} \beta_{i_{k}}$
	\EndFor
 \State {\bf Output:}  $x^N$
\end{algorithmic}
\end{algorithm}

Note that for Algorithm~\ref{algo:CD} the \gav{Assumption}~\ref{main_prop} holds with $C_n = O(n)$ (for $x \in \R^n$) and\footnote{\avg{Strictly speaking, such a constant takes place for non-adaptive variant of the CDM with specific choice of $i_k$ \cite{nesterov2012efficiency}: $\pi(i_k = j) = \frac{\beta_j}{\sum_{j'=1}^n \beta_{j'}}$. For described RACDM the analysis is more difficult \cite{pasechnyuk2021accelerated}.}} $L_{f} = \overline{L}_{f} := \frac{1}{n}\sum\limits_{i=1}^n \beta_i$ (the  average  value  of the  directional  smoothness  parameters).

As one of the motivational example, consider the following minimization problem
\begin{equation}
\min\limits_{x \in \mathbb{R}^n}~~f(x) = \gamma \ln \left( \sum_{i=1}^m \exp\left(\frac{\left[A x\right]_i}{\gamma}\right) \right) - \langle b, x \rangle,
\end{equation}
where $A \in \mathbb{R}^{m \times n}$, $b \in \mathbb{R}^n$. We denote the $i^\text{th}$ row of the matrix $A$ by $A_i$. $A$ is sparse, i.e. average number of nonzero elements in $A_i$ is less than $s$. $f$ is $L_f$-smooth w.r.t. $\|\cdot\|_{2}$ with $L_f = \max_{i=1,...,m} \|A_i\|_2^2$ and its gradient is component-wise $\beta_j$-continuous with $\beta_j = \max_{i=1,...,m} \abs{A_{i j}}$.

Fast Gradient Method (FGM)~\cite{nesterov2018primal} requires $\displaystyle O\left(\sqrt{\tfrac{L_f R^2}{\varepsilon}}\right)$
iterations with the complexity of each iteration
$\displaystyle O\left(ns\right)$. 

Coordinate Descent Method (CDM)~\cite{bubeck2015convex}
requires $\displaystyle O\left(n\tfrac{\overline{L}_f R^2}{\varepsilon}\right)$
iterations with the complexity of each iteration\footnote{Here one should use a following trick in recalculation of $\ln\left(\sum_{i=1}^m \exp \left([ A x]_i\right)\right)$ and its gradient (partial derivative). From the structure of the method we know that $x^{new} = x^{old} + \delta e_i$, where $e_i$ is $i$-th orth. So if we've already calculate $A x^{old}$ then to recalculate $A x^{new} = A x^{old} + \delta A_i$ requires only $O(s)$ additional operations independently of $n$ and $m$.}
$\displaystyle O\left(s\right)$.

Accelerated Coordinate Descent Method (ACDM)~\cite{nesterov2017efficiency,gasnikov2016accelerated}
requires $\displaystyle O\left(n\sqrt{\tfrac{\widetilde{L}_f R^2}{\varepsilon}}\right)$
iterations with the complexity of each iteration
$\displaystyle O\left(n\right)$, where $\displaystyle \widetilde{L}_f = \frac{1}{n} \sum_{j=1}^n \sqrt{\beta_j}$.

For proposed in this paper approach we have
$\displaystyle O\left(n\sqrt{\tfrac{\overline{L}_fR^2}{\varepsilon}}\right)$
iterations of CGM with complexity of each inner iteration
$\displaystyle O(s)$ and complexity of each outer iteration $\displaystyle O(ns)$. However, outer iteration executes ones per $\sim n$ inner iterations, so average-case iteration complexity is $O(s)$.

We combine all these results in the table below. From the table one can conclude that if $\overline{L}_f < L_f$, then our approach has better theoretical complexity.

\begin{table}[ht]
\centering
\begin{tabular}{|c|c|c|ll}
\cline{1-3}
 Algorithm  & Complexity&  Reference&   &  \\ \cline{1-3}
FGM & $ O\left(ns\sqrt{\frac{L_f R^2}{\varepsilon}}\right)$ & \cite{nesterov2018primal} &  &  \\ \cline{1-3}
CDM & $ O\left(ns\frac{\overline{L}_f R^2}{\varepsilon}\right)$ & \avg{\cite{nesterov2012efficiency,bubeck2015convex}} &  &  \\ \cline{1-3}
ACDM & $ O\left(n^2\sqrt{\frac{\widetilde{L}_f R^2}{\varepsilon}}\right) $ & \cite{nesterov2017efficiency} &  &  \\ \cline{1-3}
Catalyst CDM & $O\left(ns\sqrt{\frac{\overline{L}_f R^2}{\varepsilon}}\right)$ & this paper &  &  \\ \cline{1-3}
\end{tabular}
\end{table}

Note that the use of Component Descent Method allows us to improve convergence estimate by factor $\sqrt{n}$ compared to Fast Gradient Method. Indeed, for this problem we have $L_f = \max_{i=1,...,m} \|A_i\|^2_2 = O(n),$
and on the other hand $\overline{L}_f = \frac{1}{n} \sum_{j=1,...,n} \max_{i=1,...,m} \abs{A_{i j}} = O(1)$.
Therefore, the total convergence estimate for Fast Gradient Method can be written as
$$O\left(n s \cdot \sqrt{n} \cdot \sqrt{\frac{R^2}{\varepsilon}}\right),$$
and for proposed in this paper method factor $\sqrt{n}$ is reduced to $O(1)$ and could be omitted:
$$O\left(n s \cdot \sqrt{\frac{R^2}{\varepsilon}}\right).$$

The best complexity improvement is achieved if $L_f = n$, which means there is at least one row in the matrix such that $A_i = \mathbbm{1}^n$, even though all other rows can be arbitrary sparse.


\subsection{Alternating Minimization}
Consider the following problem
\begin{equation*}
\min_{x\gav{=(x_1,...,x_p)^T\in \otimes_{i=1}^p \R^{n_i}}}
f(x), 
\end{equation*}
where $f(x)$ is a $L_f$--smooth convex function (its gradient is Lipschitz continuous w.r.t. $\|\cdot\|_{2}$ with the constant $L_{f}$).

For the general case of number of blocks $p \geqslant 2$ the Alternating Minimization algorithm may be written as Algorithm \ref{AM}. There are multiple common block selection rules, such as the cyclic rule or the Gauss--Southwell rule \citep{karimi2016linear,beck2017first,diakonikolas2018alternating,tupitsa2019alternating}.

\begin{algorithm}[ht]
\caption{Alternating Minimization}
\label{AM}
\begin{algorithmic}
   \State {\bf Parameters:} Starting point $x^0$.
   \For{$k = 0,1,\ldots,N-1$}
	\State Choose $i \in \{1, \ldots, p\}$
	\State Set $x^{k+1}=\argmin\limits_{x_i}
f(\gav{x_1^k,...x_{i-1}^k,x_i,x_{i+1}^k,...,x_p^k})$

    \EndFor
 \State {\bf Output:} $x^N$
\end{algorithmic}
\end{algorithm}

Note that for Algorithm~\ref{AM} the \gav{Assumption}~\ref{main_prop} holds with $C_n = O(p)$ ($p$ -- number  of  blocks) and $L_{f}$ is the Lipschitz constant of the gradient of function $f$.

\subsection{Local Stochastic  Gradient Descent}
Local SGD \citep{khaled2019better} becomes popular first of all due to the application in federated learning \citep{kairouz2019advances}. In the core of the approach lies a very simple idea \citep{nesterov2008confidence,dvurechensky2018parallel,stich2018local,khaled2019better}: to solve considered stochastic convex optimization problem in parallel on $M$ nodes via Adaptive Stochastic Gradient Descent (SGD) \citep{ogaltsov2019adaptive} with synchronization after each $\tau$ iterations of SGD and sharing an average point. The large we want to choose $M$ the smaller we should choose $\tau$ to conserve the total (optimal) number of oracle $T$ (stochastic gradient) calls (on $M$ nodes). Say, for strongly convex case $M\tau\simeq T$ \citep{khaled2019better}. It is well known that for stochastic convex optimization problems from the oracle complexity point of view there is no difference between accelerated schemes and non-accelerated ones \citep{nemirovsky1983problem}. On the other hand, if we reduce the variance by batching acceleration could play a significant role \citep{woodworth2018graph,gorbunov2019optimal}. That is in parallelized architecture the accelerated schemes are dominant. So, the natural question: Can we accelerate local SGD? Below we'll try to demonstrate the numerical possibility of acceleration local SGD by proposed M-S Catalyst scheme. \gav{From the theoretical perspective\dmiv{,} an acceleration of local SGD is still an open problem \cite{woodworth2020local} rather than acceleration of ordinary SGD (see, for example, \cite{GasnikovBook} and reference \dmiv{therein}).}

\begin{algorithm}[ht]
\caption{Local SGD algorithm}
\label{alg:lsgd_alg}
\begin{algorithmic}[ht]
	\State{\bf Parameters:}
	$x^0 \in \mathbb{R}^n, w$~--- number of workers, $L, \mu$,
	
	$\mathcal{S}_N$~--- set of synchronization steps indices
	
	$\tau$~--- maximum difference between two consequent integers in $\mathcal{S}_N$

	\State Initialize variables $x^0_h = x^0$ for $h \in [1, w]$
	\For{$k = 0,1, \ldots, N-1$}
        \For{$h \in \{1, \dots, w\}$}~\textbf{in parallel}
        \State Sample $i^k_h$ uniformly in $[1, m]$
        \If {$k + 1 \in \mathcal{S}_N$}
        \State $x^{k+1}_h = \frac{1}{w} \sum_{j = 1}^w \left(x^k_j - \eta_k \nabla f_{i^k_h}(x^k_j)\right)$
        \Else
        \State $x^{k+1}_h = x^k_h - \eta_k \nabla f_{i^k_h}(x^k_h)$
        \EndIf
        \EndFor
	
    \EndFor

\State {\bf Output:} $\hat{x}^N=\frac{1}{w S_N} \sum_{h=1}^w \sum_{k=0}^{N-1} \xi^k x^k_h$, where $\xi^k = (\max\{16L/\mu, \tau\} + k)^2$, $S_N = \sum_{k=0}^{N-1} \xi^k$.
\end{algorithmic}	
\end{algorithm}

\subsection{Theoretical Guarantees}

Let us present the table that establishes the comparison of rates of  convergence for the above algorithms before and after acceleration via Algorithm~\ref{alg1}. In non-accelerated case these algorithms apply to the convex but non-strongly convex problem, therefore, we use estimates for the convex case from \gav{Assumption}~\ref{main_prop}. But in the case of acceleration of these methods, we apply them to a regularized function which is strongly convex. Denote $\avg{\chi} = \max \left(\sqrt{\frac{L_u}{L_{f}}}, \sqrt{\frac{ L_{f}}{L_d}}\right)$, then we represent the following table.

\begin{table}[ht]
\centering
\title{Theoretical Guarantees with and without M-S acceleration}
\begin{tabular}{|c|c|c|ll}
\cline{1-3}
 & non-accelerated & M-S  accelerated  &  &  \\ \cline{1-3}
Steepest Descent & $ \frac{L_{f}R^2}{\varepsilon}$ & $\avg{\chi} \sqrt{\frac{L_{f}R^2}{\varepsilon}}$ &  &  \\ \cline{1-3}
Random Adaptive Coordinate Descent Method & $  n \cdot \frac{\overline{L}_{f} R^2}{\varepsilon}$ & $n \cdot \avg{\chi} \sqrt{\frac{ \overline{L}_{f}R^2}{\varepsilon}}$ &  &  \\ \cline{1-3}
Alternating Minimization & $ p \cdot \frac{L_{f}R^2}{\varepsilon}$ & $p \cdot \avg{\chi} \sqrt{\frac{L_{f}R^2}{\varepsilon}}$ &  &  \\ \cline{1-3}
\cline{1-3}
\end{tabular}
\end{table}

\section{Experiments}\label{4}
In this section, we perform experiments for justifying the acceleration of the aforementioned methods in practice.
For all figures below, the vertical axis measures functional suboptimality \text{$f(x) - f(x_{\star})$ in the logarithmic scale.}
\subsection{Steepest Descent Acceleration}
\begin{figure}[ht]
\centering
     \subfloat[][Logistic regression~\eqref{eq:logistic} with rcv1\_train dataset.]{\includegraphics[width = 0.47\linewidth]{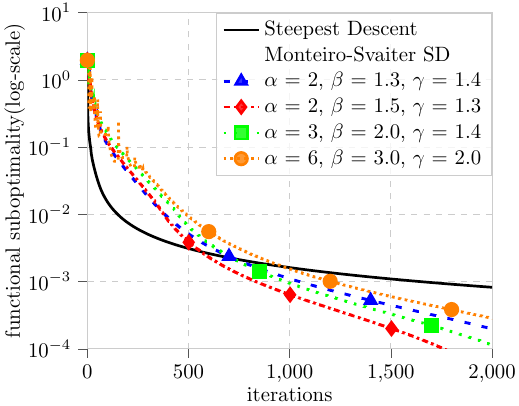}\label{fig:rcv1}}
     \hfill
     \subfloat[][Logistic regression~\eqref{eq:logistic} with a1a dataset.]{\includegraphics[width = 0.47\linewidth]{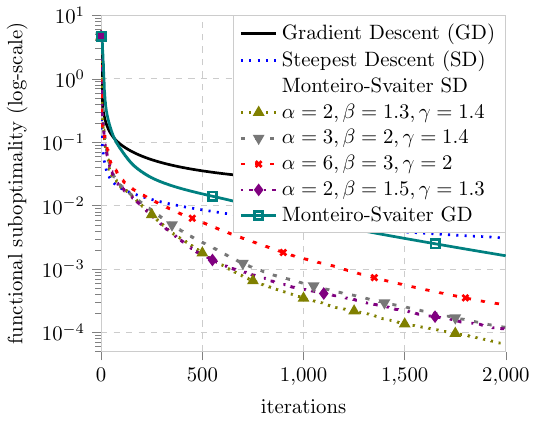}\label{fig:a1a}}
\caption{Results of Steepest Descent acceleration for different problems}
\end{figure}
In this experimental setup we consider logistic loss minimization problem
\begin{equation}\label{eq:logistic}
\min_{x\in \mathbb{R}^n}~~f(x) = \frac{1}{m}\!\sum_{j=1}^m \log(1 \!+\! \exp(-y_j z_j^{\top} x))
\end{equation}
with two different datasets from LIBSVM~\citep{chang2011libsvm} repository (\emph{rcv1\_train} and \emph{a1a}) We selected logistic regression as the objective as far as logistic regression converges with sub-linear rate like general non-strongly function that is important assumption for accelerated versions of algorithm to be provable.

In this setup we present our experimental results for {\textit{Steepest Descent}} (Algorithm~3) and its accelerated via Algorithm~2 versions with different tuples of parameters ($\alpha, \beta, \gamma$) to show the dependence of the algorithm on parameters.\footnote{For all runs with steepest descent we used $L_d = 0.0001 L_f$ and $L_u = L_f$, where $L_f$ is a real Lipschitz constant of $\nabla f$.}

In Figure~\ref{fig:rcv1} we present functional suboptimality vs aggregated amount of gradient computations (oracle calls) in logarithmic scale. To be more precise, we present the following: for every ``restart'' we split all the points into two groups. First group -- points from the inner algorithm run with ``optimal'' $L_{k+1}$. Second group is for the points, that are extra (cost of adaptation) and for this points we plot the value in point $y^k$ from the previous restart (to have the horizontal lines on plots in cases when adaptation takes so long). In the end of the day, for each restart we, first, plot ``horizontal line'' for all points from the second group and after we present points from the first group with the corresponding to them values.

As we could see from the plot, acceleration happens when M-S
acceleration scheme is used together with steepest descent but is highly dependent on the parameters of Monteiro--Svaiter algorithm. For instance, big $\alpha$ and $\beta$ make it harder to algorithm to adapt to the current ``optimal'' value of $L_{k+1}$ that makes algorithm slower. Second, selecting big $\gamma$ is not reasonable too as far as it corresponds to the big fluctuation of $L_{k+1}$ during every restart. 
Moreover, selecting $\alpha$ and $\beta$ close to each other also tends to slow down the convergence process. 

\mitya{Let us give some intuition why this happens. Let us recall, that parameters $\alpha$ and $\beta$ impact on the speed of adaptation. More precisely, the decrease of estimation $L_{k+1}$ after the one iterate is by factor of $\beta^{p}/\alpha$, where $p\in \mathbb{Z}_{+}$. This implies two different things: 
\begin{itemize}
    \item[] if $L_{k+1}< L_f$ the ``optimistic'' amount of adaptation rounds is $\log_{\alpha}{\frac{L_f}{L_{k+1}}}$, that is very big if $\alpha$ is close to $1$;
    \item[] if $L_{f} \in \left(\frac{L_{k+1}}{{\beta}}, L_{k+1}\right)$ the ``worst'' amount of adaptation rounds is $\log_{\frac{\beta}{\alpha}}{{\beta}}$, that is very big if $\alpha$ is close to $\beta$.\footnote{If $L_f < L_{k+1}/\beta $ then the scaling by factor of $\beta$ give us the situation described  in this case.}
\end{itemize}
Combining these one can find the explanation of the dependence between $\alpha$, $\beta$, and the rate of convergence.}

In Figure \ref{fig:a1a}, we add also Gradient Descent algorithm to the list of presented algorithms. To be precise enough, we use Monteiro--Svaiter acceleration without any adaptation for $L_{k+1}$. It means, that fixed constant $L_{k+1} = L_f$ is used during algorithm run.

As we could see from the set of hyper parameters $(\alpha,\beta,\gamma) = (6,3,2)$ again leads to the slowest version of accelerated algorithm. All the other set ups, also have roughly the same behavior.
Finally, for Gradient Descent acceleration also takes place and even makes it faster than Steepest Descent without acceleration. An important thing to notice is the following: Monteiro--Svaiter acceleration for adaptive algorithms (Steepest Descent) makes them faster than acceleration  of non-adaptive (Gradient Descent) in spite of cost for adaptation.

\subsection{RACDM Acceleration}

\begin{wrapfigure}[17]{R}{0.47\linewidth}
    \vspace*{-15pt}
    \includegraphics[width = \linewidth]{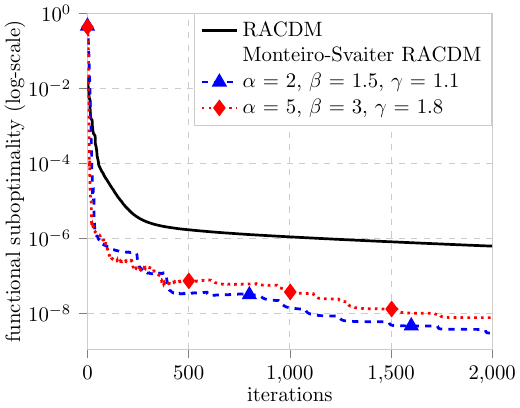}
    \caption{Results of RACDM acceleration for quadratic problem~\eqref{eq:quadratic} with Hilbert matrix, $n = 1000$.} 
    \label{fig:hilbert_plot}
\end{wrapfigure}
Let us consider quadratic optimization problem
\begin{equation}\label{eq:quadratic}
\min\limits_{x \in \mathbb{R}^n}~~f(x) = \frac{1}{2} x^{\top} A x ,
\end{equation}
for Hilbert matrix~\citep{hilbert1894} with such entries $A_{ij} = \frac{1}{i + j - 1}$. This is an example of a Hankel matrix and is known to be very ill-conditioned (e.g. for $n=6$ condition number $\approx 1.5 \cdot 10^7$~\citep{polyak1987introduction}). It leads to a degenerate optimization problem, typically very hard for gradient methods.

In Figure~\ref{fig:hilbert_plot} we compare the performance of the method~\ref{algo:CD} and its M-S accelerated version with different sets of parameters $(\alpha, \beta, \gamma)$ for problem~\eqref{eq:quadratic}.

For the horizontal axis we use number of partial derivative evaluations divided by dimensionality $n$ of the problem.
Our warm start strategy includes running inner method from the last point $y_k$ and with estimates $\hat{\beta_i}$ of smoothness constants obtained from the previous outer (M-S) iteration. The initial points $y_0=z_0$ entries
were sampled from the standard uniform distribution $\mathcal{U}(0,1)$. $L_0$ was initialized as $0.5 L_f$ and $L_d = 0.001 L_f, L_u = 100 L_f, \beta^0_i = 1/L_0$. \mitya{Here we \gav{introduce}
the relationship between 
\gav{$L_d, L_u$}\dmiv{,}
and $L_f$ only from the theoretical interest; however, the dependence between $L_d, L_u$, and $L_f$ is never used in the algorithm.}\\ 

Consider a simple case of quadratic optimization problem \eqref{eq:quadratic} with matrix $A = S^\top D S$ such that $S$ is a random orthogonal matrix. The elements of diagonal matrix $D$ are sampled from standard uniform distribution $\mathcal{U}(0, 1)$ and one random $D_{ii}$ is assigned to zero to guarantee that the smallest eigenvalue of the resulting matrix $A$ is smaller than $10^{-15}$ and thus the optimization problem is convex but not strongly-convex (up to machine precision).\footnote{Frankly speaking, for such objective functions we observe that non-accelerated  gradient descent based algorithms converge with linear rate, because of the quadratic nature of the problem and specificity of spectrum. Since $n = 100$ in these experiments we typically have that the next eigenvalue after zero is about $0.01$. This value determined the real rate of convergence} 

\begin{figure}[t]
\centering
     \subfloat[][Synthetic quadratic problem~\eqref{eq:quadratic} with matrix $A, n=100$.]{\includegraphics[width = 0.47\linewidth]{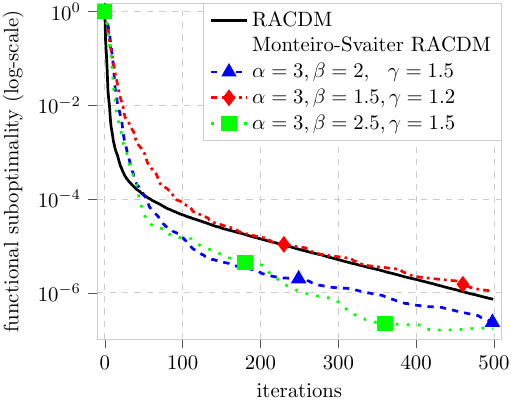}\label{fig:synth_quadratic}}
     \hfill
     \subfloat[][Logistic regression problem~\eqref{eq:logistic2} for \textit{madelon} dataset.]{\includegraphics[width = 0.47\linewidth]{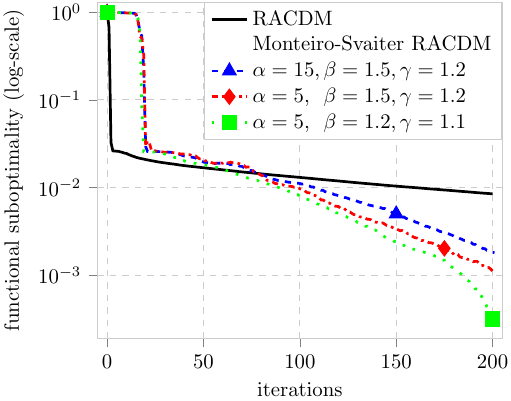}\label{fig:RACDM_madelon}}
\caption{Results of RACDM acceleration for different problems}
\end{figure}

In Figure~\ref{fig:synth_quadratic} we compare the performance of RACDM 
and its M-S accelerated version with different sets of parameters $(\alpha, \beta, \gamma)$.  

For the horizontal axis we use number of partial derivative evaluations divided by dimensionality $n$ of the problem.
Our warm start strategy includes running inner method from the last point $y_k$ and with estimates $\hat{\beta_i}$ of smoothness constants obtained from the previous iteration. The initial points $y_0=z_0$
are sampled from the standard uniform distribution $\mathcal{U}(0,1)$. $L_0$ was initialized as $1.6 L_f$ and $L_d = 0.005 L_f, L_u = 10 L_f, \beta^0_i = 1/L_0$. 
We observe that clear acceleration can be achieved not for all sets of parameters. Concretely, both $\beta$ and $\gamma$ affected convergence as one can see from the plot. Besides, we find out that for higher accuracy the proposed method can show unstable performance.

Note, that we can obtain provable acceleration by the proposed Adaptive Catalyst procedure only for convex problems. For strongly convex problems, this is no longer true either in theory or in our experiments. The reason is that the proposed M-S accelerated envelop\dmiv{e} doesn't take in to account possible strong convexity. Moreover, as far as we know, this is still an open problem, to propose a fully adaptive accelerated algorithm for strongly convex problems. The problem is in the strong convexity parameter. In practice, we met this problem in different places. For example, when we choose the restart parameter for conjugate gradient methods or try to propose accelerated (fast) gradient methods that do not require any information about strong convexity parameter but know all other characteristics of the problem.\\
    
Consider logistic loss minimization problem
\begin{equation}\label{eq:logistic2}
\min_{x\in \mathbb{R}^n}~~f(x) = \frac{1}{m}\!\sum_{j=1}^m \log(1 \!+\! \exp(-y_j z_j^{\top} x)).
\end{equation}

In Figure~\ref{fig:synth_quadratic} we compare the performance of RACDM 
and its M-S accelerated version with different sets of parameters $(\alpha, \beta, \gamma)$ for logistic regression problem with \textit{madelon} dataset from LIBSVM~\citep{chang2011libsvm} repository. Initial parameters for this set up are $L_0 = L_f, L_d = 10^{-5}, L_u = 100 L_f$. Warm start strategy and depicted values for theot horizontal axis are the same as for the quadratic problem. It is important to mention that gradient norm computation for checking Monteiro--Svaiter condition $\left\|\nabla F_{L, x}^{k+1}\left(y^{k+1}\right)\right\|_{2} \leq \frac{L_{k+1}}{2}\left\|y^{k+1}-x^{k+1}\right\|_{2}$ and full gradient step from the outer loop $z^{k+1}=z^{k}-a_{k+1} \nabla f\left(y^{k+1}\right)$ (according to Algorithm 2) were counted as evaluation of $n$ partial derivative.
For this case, we also noted that $L_0$ initialization affects the convergence significantly at the beginning and it has to be chosen lower than in the previous cases.\\

Consider also the following softmax function minimization problem
\begin{equation}\label{eq:softmax}
\min\limits_{x \in \mathbb{R}^n}~~f(x) = \gamma \ln \left( \sum_{i=1}^m \exp\left(\frac{\left[A x\right]_i}{\gamma}\right) \right) -\langle b,x \rangle,
\end{equation}
where matrix $A \in \mathbb{R}^{m \times n}$ is such that average number of nonzero elements in $A_i$ is less than $s \ll m$ and one of the rows $A_k$ is non-sparse. $f$ is Lipschitz smooth with the constant $L_f = \max_{i=1,...,m} \|A_i\|_2^2$ and has component-wise Lipschitz continuous gradient with constants $\beta_j = \max_{i=1,...,m} \abs{A_{i j}}$.

\begin{figure}[!ht]
\begin{minipage}{0.47\textwidth}
    \centering
    \includegraphics[width=\linewidth]{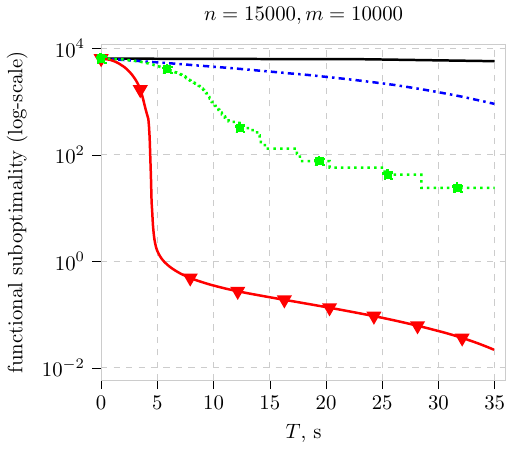}
    \caption{Softmax problem~\eqref{eq:softmax} with random sparse matrix.}
    \label{fig:umc1}
\end{minipage}
\hfill
\begin{minipage}{0.47\textwidth}
    \centering
    \includegraphics[width=\linewidth]{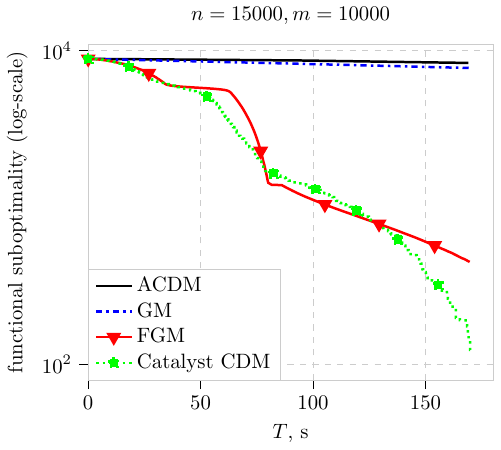}
    \caption{Softmax problem~\eqref{eq:softmax} with heterogeneous sparse matrix.}
    \label{fig:umc2}
\end{minipage}
\end{figure}

In figures \ref{fig:umc1} and \ref{fig:umc2} we compare the performance of the Accelerated Coordinate Descent Method (ACDM) \avg{\cite{nesterov2017efficiency}}, Gradient Method (GM) \avg{\cite{nesterov2018lectures}}, Fast Gradient Method (FGM) \avg{\cite{nesterov2018lectures}}, and proposed approach -- accelerated via Algorithm~\ref{alg1} variant of \avg{non-adaptive} Coordinate Descent Method (Catalyst CDM) for strongly convex w.r.t. $\|\cdot\|_2$ auxiliary problem, in which $i_k$ is drawn from the distribution $\pi$ defined by \avg{\cite{nesterov2012efficiency}} 
$$
\pi(i_k = j) = \frac{\beta_j}{\sum_{\avg{j'}=1}^n \beta_{\avg{j'}}}.
$$

In the case of randomly-generated $A$ with $A_{i j} \in \{0, 1\}$, $s \approx 0.2 m$, non-sparse row $A_k$ with uniformly random components and $\gamma=0.6$ proposed method converges faster (in terms of working time) than all comparable methods except FGM. However, in other setting, for the smaller count of nonzero elements in the matrix ($s \approx 0.75 m$), non-sparse row $A_k = \mathbbm{1}^n$, and with higher variation in the sparsity of the rows of $A$ ($0.9 m$ rows with $0.1 n$ nonzero elements and $0.1 m$ rows with $0.9 n$ nonzero elements) --- proposed method converges faster than FGM due to $\overline{L}_f \ll L_f$.

\subsection{Alternating Least Squares Acceleration}

Consider the typical collaborative filtering problem: completion of the user-item preferences matrix with estimated values based on a little count of observed ratings made by other users on other items. The considered being accelerated AM algorithm is induced by the idea of preferences matrix factorization and estimating unknown rating $r_{ui}$ associated with the user $u$ and the item $i$ as a product $x_u^\top y_i$, where the vectors $x_u$ and $y_i$ are being optimized variables. Following the approach set out in \citep{hu2008collaborative}, formulate such an optimization problem:

\begin{equation}\label{eq:mcomp_problem}
\min\limits_{x, y}~~F(x, y) = \sum_{u, i} c_{ui} \left(p_{ui} - x_u^\top y_i\right)^2 +  \lambda \left(\sum_{u} ||x_u||_2^2 + \sum_{i} ||y_i||_2^2\right),
\end{equation}
where $c_{ui}$ is confidence in observing $r_{ui}$, in our case expressed as $c_{ui} = 1 + 5 r_{ui}$, $p_{ui}$ is binarized rating:
$$
p_{ui} = \left\{
                \begin{array}{ll}
                  1 \quad r_{ui} > 0,\\
                  0 \quad r_{ui} = 0,\\
                \end{array}
              \right.
$$
and $\lambda \left(\sum_{u} ||x_u||_2^2 + \sum_{i} ||y_i||_2^2\right)$~--- regularization term preventing overfitting during the learning process (in our case, the regularization coefficient is set to $\lambda = 0.1$).

For described objective functional optimization we used modified Algorithm~\ref{AM} in that on every iteration functional optimizing with relation to $x$ and $y$ consistently (for that functional we can get the explicit expression for the solutions of equations $\nabla_{x} f(x, y) = 0$ and $\nabla_{y} f(x, y) = 0$, computational effective matrix expressions for solutions of these auxiliary problems are presented in \citep{hu2008collaborative}).

The considered objective function is not convex, so instead of the described Adaptive Catalyst scheme for accelerating should be used the modified one, in which the step of updating variable $y_k$ replaced with such construction:

\begin{gather}
\tilde{y}_{k+1} \approx \argmin_{y} F^{k+1}_{L, x}(y) \\
y_{k+1} = \arg \min \left\{f(y)\;|\;y \in \{y_k, \tilde{y}_{k+1}\}\right\}
\end{gather}

Used in experiments sparse matrix $\{r_{ui}\}_{u, i}$ generated from \textit{radio} dataset\footnote{https://www.upf.edu/web/mtg/lastfm360k} with ratings made by listeners on compositions of certain artists. Count of considered users was 70, count of artists~--- 100, and sparsity coefficient of the matrix was near the $2\%$.

\begin{figure}[!ht]
   \begin{minipage}{0.47\textwidth}
    \centering
    \includegraphics[width=\linewidth]{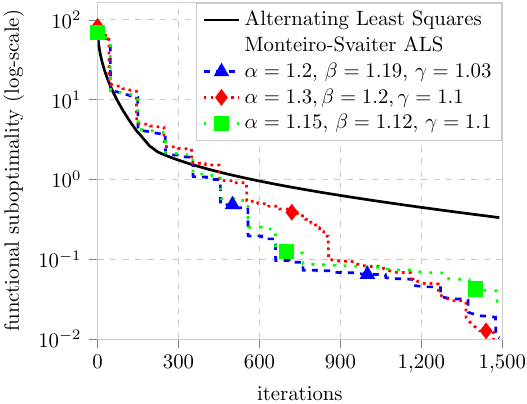}
    \caption{Matrix completion problem~\eqref{eq:mcomp_problem} with different $(\alpha, \beta, \gamma)$.}
    \label{fig:matrix_completion}
    \end{minipage}
   \hfill
   \begin{minipage}{0.47\textwidth}
     \centering
     \includegraphics[width=\linewidth]{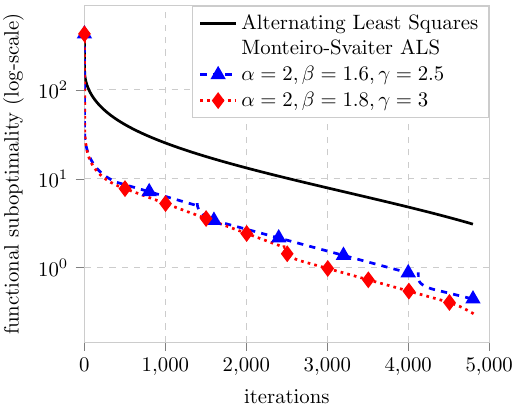}
     \caption{Matrix completion problem~\eqref{eq:mcomp_problem} with big matrix.}\label{fig:matrix_completion2}
   \end{minipage}
\end{figure}

In Figure~\ref{fig:matrix_completion} we compare the performance of the modified Algorithm~\ref{AM} and their accelerated via Algorithm~\ref{alg1} versions (with a different choice of hyperparameters $(\alpha, \beta, \gamma)$). 

The horizontal axis measures the number of variables recomputations. The plot show that there was the acceleration of the base algorithm and both parameters $\beta$ and $\gamma$ had an impact on the convergence rate.

In addition, consider the performance of the Alternating Least Squares algorithm for the problem with a larger being estimated matrix $\{r_{ui}\}_{u, i}$.

In Figure~\ref{fig:matrix_completion2} we compare the performance of the Alternating Least Squares algorithm and its accelerated via Monteiro--Svaiter algorithm versions for problem \ref{eq:mcomp_problem} with $\lambda = 5$ and matrix $\{r_{ui}\}_{u, i}$ of the size 150 users $\times$ 300 artists, generated from radio dataset.
The horizontal axis measures the number of variables recomputations. The plot shows that there was the acceleration of the base algorithm and both parameters $\beta$ and $\gamma$ had an impact on the convergence rate.

\subsection{Local SGD Acceleration}
\begin{figure}[h]
   \begin{minipage}{0.47\textwidth}
    \centering
    \includegraphics[width=\linewidth]{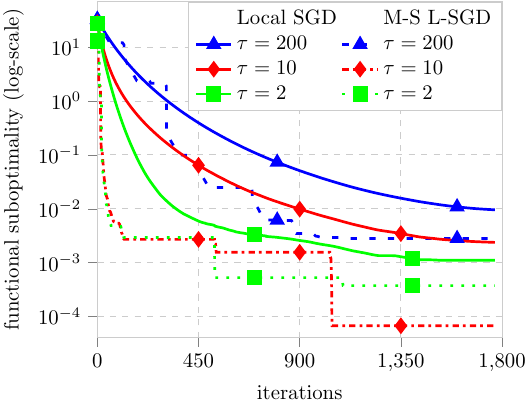}
    \caption{Regularized logistic loss~\eqref{eq:logloss} for different synchronization intervals $\tau$.}
    \label{fig:lsgd_1}
    \end{minipage}
   \hfill
   \begin{minipage}{0.47\textwidth}
     \centering
     \includegraphics[width=\linewidth]{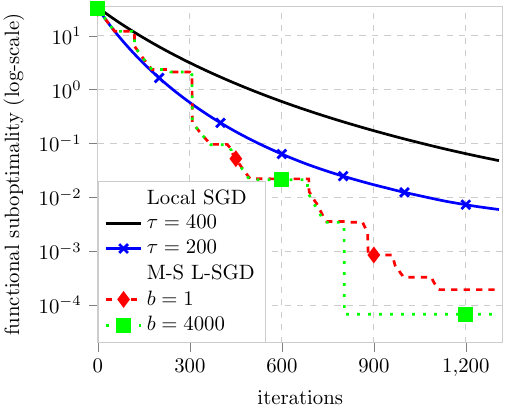}
     \caption{Logistic loss~\eqref{eq:logloss} minimization with minibatching.}\label{fig:lsgd_2}
   \end{minipage}
\end{figure}

Consider the following \dmiv{$\ell_2$}-regularized logistic loss minimization problem:
\begin{equation}\label{eq:logloss}
\min_{x\in \mathbb{R}^n}~~F(x) = \frac{1}{m} \sum_{j=1}^m \log(1 + \exp{(-y_j x^\top p_j)}) + g(x)
\end{equation}

where the features matrix $P \in \mathbb{R}^{m\times n}$, labels $y \in \{0, 1\}^m$
and $g$ is a regularization term, aggregated by two different regularizers for the sparse features $I_S \subset [1, n]$ (with coefficient $\lambda_1$) and the dense features $I_D \subset [1, n], I_S \cap I_D = \varnothing$ (with coefficient $\lambda_2$):
$$
g(x) = \lambda_1 \sum_{i \in I_S} x^2_i + \lambda_2 \sum_{i \in I_D} x^2_i.
$$
In this experiment, we use the \textit{adults} dataset with one-hot encoded categori\dmiv{c}al features and binarized `work class` feature as a label, $m = 40000$, $n=100$, $\lambda_1 = 1.1$, $\lambda_2 = 2.1$, the sparsity coefficient (percentage of features with a fraction of zeros less than $0.2$) is equal to $4\%$. 
Also, the initial value of the model's weights\dmiv{  is} randomly generated from the uniform distribution $x_{0, i} \sim \mathcal{U}(0, 1)\;\forall i \in [1, n]$. 

In Figure~\ref{fig:lsgd_1} we compare the performance of the Local SGD and its accelerated via Algorithm~\ref{alg1} versions. Parameters used: $w = 20$ (amount of computing nodes), varying $\tau$ (number of iterations between two consequent synchronization steps), $\alpha=1.15, \beta=1.12, \gamma=1.1$ (for Monteiro--Svaiter).

The horizontal axis measures the number of outer iterations (one outer iteration includes recomputation of the variables in all computing nodes). The plot shows that there was the acceleration of the base algorithm and synchronization interval had an impact on the convergence rate.

In Figure~\ref{fig:lsgd_2} we compare the performance of the Algorithm~\ref{alg:lsgd_alg} and its Monteiro--Svaiter (with parameters $\alpha=1.15, \beta=1.12, \gamma=1.1$ and $\tau=200$) accelerated versions ($w = 20$ and $\tau \in \{200, 400\}$) for problem~\eqref{eq:logloss} with applying the minibatch technique \dmiv{(where parameter $b$ controls the batch size)}. 
The horizontal axis measures the number of outer iterations (one outer iteration includes the recomputation of the variables in all the computing nodes). The plot shows that there was the acceleration of the base algorithm and, moreover, that using a batch of samples instead of one sample for calculating the stochastic gradient can improve the convergence rate.

\section*{Conclusion}
In this work\dmiv{,} we present the universal framework for accelerating the non-accelerated adaptive methods such as Steepest Descent, Alternating Least Squares Minimization\dmiv{,} and RACDM and show that acceleration works in practice (code is available online on \href{https://github.com/AdaptiveCatalyst/Adaptive-Catalyst}{GitHub}). Moreover, we show theoretically that for \dmiv{the }non\dmiv{-}adaptive run proposed in this paper\dmiv{,} acceleration has \gav{in a log-factor} better rate than via Catalyst. \gav{Note, that this ``fight'' for the log-factor in accelerated procedure's become popular in the last time, see \cite{kovalev2020optimal,li2020optimal} for concrete examples. In this paper\dmiv{,} we eliminate log-factor in a rather big generality.}
\section*{Acknowledgements}
We would like to thank Soomin Lee (Yahoo), Erik Ordentlich (Yahoo), César A. Uribe (MIT), Pavel Dvurechensky (WIAS, Berlin) and Peter Richtarik (KAUST) for useful remarks. \gav{We also would like to thank \egor{anonymous} reviewers for their fruitful comments.}
\bibliography{bibliography.bib}
\bibliographystyle{tfs.bst}

\end{document}